\documentclass[12pt]{article}
\usepackage{amssymb}
\usepackage[all,arc]{xy}
\usepackage{mathrsfs}
\usepackage{amsfonts}
\headsep =
0.0in \headheight = 0.0in \topmargin = 0.3in %\oddsidemargin=0.1in

\usepackage{amsmath,amsthm,amsfonts,amssymb,amscd}

\allowdisplaybreaks[4]
\newtheorem{defn}{Definition}[section]
\newtheorem{them}[defn]{Theorem}
\newtheorem{lem}[defn]{Lemma}

\newtheorem{cor}[defn]{Corollary}

\newtheorem{prop}[defn]{Proposition}

\newtheorem{con}[defn]{Conjecture}

\numberwithin{equation}{section}
\newenvironment {Proof} {\noindent {\bf Proof.}}{\quad $\square$\par\vspace{3mm}}

\begin{document}
\title{On a conjecture for the signless Laplacian spectral radius of cacti with given matching
number\footnote{L. You's research is supported by the National Natural Science Foundation of China (Grant No. 11571123)
and the Guangdong Provincial Natural Science Foundation %NSF of Guangdong Province
 (Grant No. 2015A030313377), S. Li's research is supported by the National Natural Science Foundation of China (Grant No. 11271149).}}

\author{Yun Shen$^{a,}$
\footnote{{\it{Email address:\;}}1192158571@qq.com. } %; {\it{phone number:\;}}+8615626456721.}
 \ \ \ Lihua You$^{a,}$\footnote{{\it{Corresponding author:\;}}ylhua@scnu.edu.cn.} %; {\it{phone number:\;}}+8618602062229.}
  \ \ \ Minjie Zhang$^{b,}$\footnote{{\it{Email address:\;}}zmj1982@21cn.com}
 \ \ \ Shuchao Li$^{b,}$\footnote{{\it{Email address:\;}}lscmath@mail.ccnu.edu.cn}}\vskip.2cm
\date{{\small
$^{a}$ School of Mathematical Sciences, South China Normal University,\\ Guangzhou, 510631, P.R. China\\
$^{b}$ Faculty of  Mathematics and Statistics, Central China  Normal University,\\ Wuhan, 430079, P.R. China }}
\maketitle
\noindent {\bf Abstract}
A connected graph $G$ is a cactus if any two of its cycles have at most one common vertex.
Let  $\ell_n^m$ be the set of  cacti on $n$ vertices with matching number $m.$
 S.C. Li and M.J. Zhang  determined the unique graph with the maximum signless Laplacian spectral radius
among all cacti in $\ell_n^m$ with $n=2m$.
In this paper, we characterize the case $n\geq 2m+1$.
  This confirms the conjecture of Li and Zhang
  (S.C. Li, M.J. Zhang, On the signless Laplacian index of cacti with a given number of pendant vetices, Linear Algebra Appl. 436, 2012, 4400--4411).
Further, we  characterize the unique graph with the maximum signless Laplacian spectral radius
among all cacti on $n$ vertices.



\section{Introduction}

\hskip.5cm
Spectral graph theory (for example, {\cite{1997DC, 1995JMB, 2010}} et al) studies properties of graphs using the spectrum of related matrices.
In this paper, we consider only simple graphs (i.e., finite,undirected graphs without loops or multiple edges),
and follow  \cite{1979, 1976EPC, 1997DC, 1995JMB, 2010} for terminology and notations.

Let $G=(V,E)$ be a simple graph with vertex set $V=V(G)=\{v_1, v_2,\ldots,$ $v_n\}$ and edge set $E=E(G)$.
Let $A(G)=(a_{ij})$ denote the adjacency matrix of $G$, where $a_{ij}=1$ if vertex $v_i$ and vertex $v_j$  are adjacent in $G$ and $0$ otherwise.
 %The spectral radius of $A(G)$, denoted by $\rho(G)$, is called the spectral radius of $G$.
 Let $diag(G)=diag(d_1, d_2, \ldots, d_n)$ be the diagonal matrix with degree of the vertices of $G$
 and $Q(G)=diag(G)+A(G)$ be the signless Laplacian matrix of $G$. %, $L(G)=diag(G)-A(G)$ be the Laplacian matrix of $G$.
 It is well known that $A(G)$ is a real symmetric matrix and $Q(G)$ is a positive semidefinite matrix. The eigenvalues of $Q(G)$ can be ordered as
 $$q_1(G)\geq q_2(G)\geq \cdots \geq q_n(G)\geq 0,$$  where $q_1(G)$ is the largest signless Laplacian eigenvalue of $G$ and it
 is called the signless Laplacian spectral radius of $G$, denoted by $q(G)$.
 It is easy to see that if $G$ is connected, then $A(G)$ and thus $Q(G)$ is a nonnegative irreducible matrix.
 By the Perron-Frobenius theory, $q(G)$ has multiplicity one and there exists a unique positive unit eigenvector, say $\textbf{x}=(x_1,x_2,\ldots,x_n)^T,$ corresponding to $q(G),$ which is called the Perron vector of $Q(G)$.
 It will be convenient to associate a labelling of vertices of $G$ (with respect to $\textbf{x}$) in which $x_v$ is a label of $v.$
 The signless Laplacian characteristic polymomial  of $G,$ is equal to det$(xI_n-Q(G)),$   denoted by $\psi(G,x)$ ( or, for short, by $\psi(G)$).

In order to state our results, we introduce some notation and terminology.
Let $P_n,S_n$ and $C_n$ be the path, star and cycle on $n$ vertices, respectively.
Let $G-v,$ $G-uv$ denote the graph obtained from $G$ by deleting vertex $v\in V(G),$ or an edge $uv\in E(G),$ respectively
(this notation is naturally extended if more than one vertex or edge is deleted). Similarly,
Let $G+v$ be  obtained from $G$ by adding vertex $v\notin V(G)$
(note, if a vertex $v$ is added to $G,$ then its neighbours in $G$ should be specified somehow),
$G+uv$ be obtained from $G$ by adding an edge $uv \notin E(G),$ where $u,v\in V(G).$
%For $uv\in E(G),$ let $G_{uv}$ be the graph obtained from $G$ by subdividing the edge $uv,$ that is,
%by replacing $uv$ with edges $uw$ and $wv,$ where $w$ is an additional vertex.
For $v\in V(G),$ let $N_G(v)$ (or $N(v)$ for short) denote the set of all the adjacent vertices of $v$ in $G.$
%An $internal~path$ is a path or a cycle, in which the initial and terminal vertices have degree at least three
%and the internal vertices have degree two.

 The degree of a vertex $v$ in $G$ is denoted by $d_G(v).$ If $d_G(v)=1,$ then $v$ is a pendant vertex of $G$.
 %A $quasi-pendant~vertex$ is a vertex adjacent to a pendant vertex.
 An edge  associated with a pendant vertex is a pendant edge.
 Two distinct edges in a graph $G$ are independent if they do not have a common end vertex in $G.$
 A set of pairwise independent edges of $G$ is called a matching of $G,$
 while a matching of maximum cardinality is a maximum matching of $G.$
  The matching number $m$ of $G$ is the cardinality of a maximum matching of $G.$
 Let $M$ be a matching of $G.$ The vertex $v$ in $G$ is $M$-saturated if $v$ is incident with an edge in $M;$
 otherwise, $v$ is $M$-unsaturated. A perfect matching $M$ of $G$ means that each vertex of $G$ is $M$-saturated.
 Clearly, every perfect matching is maximum.

 We call graph $G$ a cactus if $G$ is connected and any two of its cycles have  at most one common vertex.
For a cactus graph $G,$ we call it a $bundle$ if all cycles of $G$ have exactly one common vertex.
Let $\ell_n^m$ denote  the set of cacti with $n$ vertices and matching number $m.$

Recently there is a lot of work on the spectral radius or the  signless Laplacian spectral radius of graphs,
see \cite{2009AM, 2007LAA, 2009IM, 2010LAA, 2009GC, 2009PIM,  2010LAA1, 2002LAA, SP, 2012LAA, 2008CM, 2010LAA3, 2010DAM,  2015You, 2009DAM} et al.
Some investigation on graphs with prefect matching or with given matching number is an important topic in the theory of graph spectra,
%The problem concerning maximal spectral radii of graph(s) with perfect matching (resp. maximum matching)
%has attracted much attention in the literature.
see  \cite{1998, 2003, 2003LAA, 2004, 2007, 2007LMA, 2009LAA, 2012, 2014, SP, 2011D, 2012LAA, 2008J, 2014SWJ, 2015PIM} et al.
Cacti has been an interesting topic in  chemical and mathematical literature, see \cite{1990JMC, 2000AAM, 2006PIM, 1969CJM, 2014, 2011D, 2009AMO, 2012F, 2012LAA, 2008DAM, 2007MCMCC, 2010DM, 2006MC, 2011MCM, 2011LAA, 2010LMA, 2013LAA, 2004DAM} et al.

Li and Zhang \cite{2012LAA} determined the unique graph with the maximum signless Laplacian spectral radius
among all cacti in $\ell_n^m$ with $n=2m$, and gave a conjecture about the case $n\geq 2m+1$ as follows.
%Motivated by these facts, we study the same question for $\ell_n^m,$ a set of polycyclic graphs (called cacti) in
%which any two of its cycles have at most one common vertex and each cactus is connected with matching number $m.$

\begin{con}\label{con11}{\rm(\cite{2012LAA}, Conjecture 3.4)}
Let $G$ be a graph in $\ell_n^m.$ Then

{\rm (i) } If $n=2m+1,$ then $q(G)\leq\frac{5+\sqrt{4n-3}}{2}$  with equality if and only if $G \cong H_m^0.$

{\rm (ii) } If $n\geq2m+2,$ then $q(G) \leq q(H_{m-1}^{n-2m+1})$ with equality if and only if $G \cong H_{m-1}^{n-2m+1},$ where $q(H_{m-1}^{n-2m+1})$ is the largest root of the equation $x^3 - (n - 2m+7)x^2 + (3n-8m+8)x +2m -8 =0.$
\end{con}

In this paper, some useful lemmas are given in Section 2. In Section 3, we characterize the unique graph with the maximum signless Laplacian spectral radius among all cacti in $\ell_n^m$ with the case $n\geq 2m+1$. This improves and confirms Conjecture  \ref{con11} of Li and Zhang in \cite{2012LAA}. Further, we  characterize the unique graph with the maximum signless Laplacian spectral radius among all cacti on $n$ vertices.

\section{Some preliminaries}

\hskip.6cm In this section, we introduce some  lemmas %basic properties
which we need to use in the presentations and proofs of our main results in Sections 3.

\begin{lem}\label{lem21}{\rm(\cite{2005DM}, Theorem 2.1)}
Let $u$ and $v$ be distinct vertices of a connected graph $G.$ Suppose that $w_1, w_2, \ldots, w_s (s\geq1)$ are neighbors of $v$ but not $u$ and they are all different from $u.$
Let $x$=$(x_1,x_2,\ldots,x_n)^T$ be the Perron vector of $Q(G),$
 $H$  obtained from $G$ by deleting the edges $vw_i$ and adding the edges $uw_i$ for $i=1,2,\ldots,s.$
If $x_v\leq x_u,$ then $q(G)<q(H).$
\end{lem}

\begin{lem}\label{lem22}{\rm(\cite{2010}, Theorem 8.17)}
Let $G$ be a connected graph with  a non-pendant edge $e=uv$ satisfying $N(u)\cap N(v)=\emptyset$.
Let $H$ be the graph obtained from $G$ by deleting edge $uv$, identifying vertex $u$ and vertex $v$, %contracting $uv$ to a vertex $w$
and adding a new pendant edge to $u(=v)$.
Then $q(G)<q(H)$.
\end{lem}

\begin{defn}\label{defn220}{\rm(\cite{1979}, Chapter 2)} %µÚ26Ò³ 1 Introduction
Let $A=(a_{ij}), B=(b_{ij})$ be $n\times n$ matrices. If $a_{ij}\leq b_{ij}$ for all $i$ and $j$, then $A\leq B$.
If $A\leq B$ and $A\neq B$, then $A< B$. %If $a_{ij} < b_{ij}$ for all $i$ and $j$, then $A \ll B$.
\end{defn}

\begin{lem}\label{lem24}{\rm(\cite{1979}, Chapter 2)}
Let $A, B$ be $n\times n$ matrices with the spectral radius $\rho(A)$ and $\rho(B)$. If $0\leq A\leq B$, then $\rho(A) \leq \rho(B)$.
Furthermore, if $0\leq A<B$ and $B$ is irreducible, then $\rho(A) < \rho(B)$.
\end{lem}

\begin{lem}\label{lem25}{\rm(\cite{1979}, Chapter 2)}
Let $m< n$, $A, B$ be $n\times n$, $m\times m$ nonnegative matrices with the spectral radius $\rho(A)$ and $\rho(B)$, respectively.
If $B$ is a principal submatrix of $A$, then $\rho(B) \leq \rho(A)$.
Furthermore, if $A$ is irreducible, then $\rho(B) < \rho(A)$.
\end{lem}

By Lemmas \ref{lem24}--\ref{lem25} and the definitions of $Q(G)$ and $q(G)$, we have the following  result in terms of graphs.

\begin{cor}\label{cor26}
Let $G$ be a connected graph. If $H$ be a subgraph of $G$, then $q(H)\leq q(G)$.
If $H$ is a proper  subgraph of $G$, then $q(H)<q(G)$.
\end{cor}

%\begin{lem}\label{lem23}{\rm(\cite{2009GC})}
%If $G_1$ is a proper subgraph of a connected graph $G,$ then $q(G_1)<q(G).$
%\end{lem}

 Let $I_p$ be the $p\times p$ identity matrix and $J_{p,q}$ be the $p\times q$ matrix in which every entry is $1$, or simply $J_p$ if $p=q$. Let $M$ be a matrix of order $n$, $\sigma(M)$ be the spectrum of the matrix $M$.
  %$P_M(\lambda)=det(xI_n-M)$ be the characteristic polynomial of matrix $M$.

\begin{defn}\label{defn27}{\rm(\cite{2014a})}
Let $M$ be a real matrix of order $n$ described in the following block form

\begin{equation}\label{eq21}
M = \left(\begin{array}{ccc}
M_{11} & \cdots & M_{1t}\\
\vdots &  \ddots      &\vdots \\
 M_{t1}& \cdots & M_{tt}\\
\end{array}\right),
\end{equation}

\noindent where the diagonal blocks $M_{ii}$ are $n_i\times n_i$ matrices for any $i\in\{1,2,\ldots, t\}$ and $n=n_1+\ldots+n_t$.
For any $i,j\in\{1,2,\ldots, t\}$, let $b_{ij}$ denote the average row sum of $M_{ij}$, i.e. $b_{ij}$ is the sum of all entries in $M_{ij}$ divided by the number of rows. Then $B(M) = (b_{ij})$ (simply by $B$) is called the quotient matrix of $M$.  In addition, if for each pair $i, j$, $M_{ij}$ has constant row sum, then $B(M)$ is called the equitable quotient matrix of $M$.
\end{defn}

\begin{lem}\label{lem28}{\rm(\cite{2015You})}
Let $M=(m_{ij})_{n\times n}$ be defined as (\ref{eq21}),  and for any $i,j \in\{ 1,2,\ldots,t\}$, the row sum of each block $M_{ij}$ be constant.
Let $B=B(M)=(b_{ij})$ be the equitable quotient matrix of $M$, and $\lambda$ be an  eigenvalue of $B$. Then $\lambda$  is also  an eigenvalue of $M$.
\end{lem}

\begin{lem}\label{lem29}{\rm(\cite{2015You})}
Let $M$ be defined as (\ref{eq21}), and for any $i,j \in\{ 1,2,\ldots,t\}$,
 $M_{ii} = l_iJ_{n_i} + p_iI_{n_i},$ $M_{ij} = s_{ij}J_{n_i,n_j}$ for $i\not= j$, where $l_i, p_i, s_{ij}$ are real numbers,
  $B=B(M)$ be the  quotient matrix of $M$.  Then
\begin{equation}\label{eq22}
\sigma(M)=\sigma(B)\cup \{p_i^{[n_i-1]} \mid i = 1,2,\ldots,t\},
\end{equation}
 where $\lambda^{[t]}$ means that $\lambda$ is an eigenvalue with multiplicity $t$.
\end{lem}

\setlength{\unitlength} {4mm}
\begin{center}
\begin{picture}(40,14)
\put(4,3){\circle* {0.3}}   \put(6,3){\circle* {0.3}}     \put(8,3){$\cdots$}   \put(11,3){\circle* {0.3}}
\put(4,2.5){$\underbrace{\quad\quad\quad\quad\quad\quad\quad}$}
\put(7.5,1){$k$}
\put(8,7){\circle* {0.3}}    \put(9,7){$v_0$}
\put(8,7){\line(-1,-1){4}} \put(8,7){\line(-1,-2){2}}  \put(8,7){\line(3,-4){3}}
\put(2,12){\circle* {0.3}}   \put(4,12){\circle* {0.3}}  \put(6,12){\circle* {0.3}}  \put(8,12){\circle* {0.3}}
\put(9.5,11.8){$\cdots$}  \put(12,12){\circle* {0.3}}  \put(14,12){\circle* {0.3}}
\put(2,12){\line(1,0){2}}  \put(6,12){\line(1,0){2}} \put(12,12){\line(1,0){2}}
\put(2,12.5){$\overbrace{\quad\quad\quad\quad\quad\quad\quad\quad\quad\quad\quad\quad}$}
\put(8,7){\line(-6,5){6}} \put(8,7){\line(-4,5){4}} \put(8,7){\line(-2,5){2}} \put(8,7){\line(0,5){5}}
\put(8,7){\line(4,5){4}} \put(8,7){\line(6,5){6}}
\put(8,13.5){$s$}
\put(5.5,-0.5){Fig.1 $H_s^k$}

\put(24,3){\circle* {0.3}}   \put(26,3){\circle* {0.3}}     \put(28,3){$\cdots$}   \put(31,3){\circle* {0.3}}
\put(24,2.5){$\underbrace{\quad\quad\quad\quad\quad\quad\quad}$}
\put(26.5,1){$k-1$}
\put(28,7){\circle* {0.3}}    \put(26,7){$v_0$}
\put(30,7){\circle* {0.3}}   \put(32,7){\circle* {0.3}}
\put(28,7){\line(1,0){2}}   \put(30,7){\line(1,0){2}}
\put(28,7){\line(-1,-1){4}} \put(28,7){\line(-1,-2){2}}  \put(28,7){\line(3,-4){3}}
\put(22,12){\circle* {0.3}}   \put(24,12){\circle* {0.3}}  \put(26,12){\circle* {0.3}}  \put(28,12){\circle* {0.3}}
\put(29.5,11.8){$\cdots$}  \put(32,12){\circle* {0.3}}  \put(34,12){\circle* {0.3}}
\put(22,12){\line(1,0){2}}  \put(26,12){\line(1,0){2}} \put(32,12){\line(1,0){2}}
\put(22,12.5){$\overbrace{\quad\quad\quad\quad\quad\quad\quad\quad\quad\quad\quad\quad}$}
\put(28,7){\line(-6,5){6}} \put(28,7){\line(-4,5){4}} \put(28,7){\line(-2,5){2}} \put(28,7){\line(0,5){5}}
\put(28,7){\line(4,5){4}} \put(28,7){\line(6,5){6}}
\put(28,13.5){$s$}
\put(25.5,-0.5){Fig.2 $L_s^k$}
\end{picture}
\end{center}

\begin{lem}\label{lem210}
Let $n$ be positive integer and $s,k$ be nonnegative integers with $2s+k+1=n$,  $H_s^k$ be a graph on $n$ vertices as in Fig. 1. Then
\begin{equation}\label{eq23}
\psi(H_s^k)=(x-1)^{\frac{n+k-3}{2}}(x-3)^{\frac{n-k-3}{2}}[x^3-(n+3)x^2+3nx-2n+2k+2].
\end{equation}
\end{lem}

\begin{Proof}
%We note that $2s+t+1=n$ and %It is obvious that  the signless Laplacian matrix of $H_s^t$,   $Q(H_s^t)$, satisfies the following,
Clearly, we have
$$Q(H_s^k) = \left(\begin{array}{ccccccccccc}
 n-1 & 1 & 1 & 1 & 1 &\cdots & 1 & 1 & 1 & \cdots & 1\\
 1 & 2 & 1 & 0  & 0 & \cdots & 0 & 0 & 0 & \cdots & 0\\
 1 & 1 & 2 & 0  & 0 & \cdots & 0 & 0 & 0 & \cdots & 0\\
 1 & 0 & 0 & 2  & 1 & \cdots & 0 & 0 & 0 & \cdots & 0\\
 1 & 0 & 0 & 1  & 2 & \cdots & 0 & 0 & 0 & \cdots & 0\\
 \vdots & \vdots & \vdots & \vdots  & \vdots & \ddots & \vdots & \vdots & \vdots & \ddots & \vdots\\
 1 & 0 & 0 & 0  & 0 & \cdots & 2 & 1 & 0 & \cdots & 0\\
 1 & 0 & 0 & 0  & 0 & \cdots & 1 & 2 & 0 & \cdots & 0\\
 1 & 0 & 0 & 0  & 0 & \cdots & 0 & 0 & 1 & \cdots & 0\\
 \vdots & \vdots & \vdots & \vdots  & \vdots & \ddots & \vdots & \vdots & \vdots & \ddots & \vdots\\
 1 & 0 & 0 & 0  & 0 & \cdots & 0 & 0 & 0 & \cdots & 1\\
  \end{array}\right),$$
then it can be written as follows:
\begin{equation}\label{eq24}
  Q(H_s^k) = \left(\begin{array}{cccccc}
  (n-2)J_1+I_1 & J_{1,2} & J_{1,2}  & \cdots & J_{1,2} &J_{1,k}\\
 J_{2,1} & J_2+I_2 & 0 & \cdots & 0  & 0\\
 J_{2,1} & 0 & J_2+I_2  & \cdots & 0 & 0\\
 \vdots &  \vdots  & \vdots & \ddots &  \vdots & \vdots\\
 J_{2,1}& 0 & 0  & \cdots & J_2+I_2  &  0\\
 J_{k,1}& 0 & 0  & \cdots & 0  &  I_k\\
  \end{array}\right).
  \end{equation}

 Let $B_1(H_s^k)$ be  the corresponding equitable quotient matrix of $Q(H_s^k)$. Then   by (\ref{eq24}) and Lemma \ref{lem29},
 we have $\sigma(Q(H_s^k))=\sigma(B_1(H_s^k))\cup\{1^{[\frac{n+k-3}{2}]}\}$, where
 $$B_1(H_s^k) = \left(\begin{array}{cccccc}
n-1 & 2 & 2  & \cdots & 2 &k\\
 1 & 3& 0 & \cdots & 0  & 0\\
 1 & 0 & 3  & \cdots & 0 & 0\\
 \vdots &  \vdots  & \vdots & \ddots &  \vdots & \vdots\\
 1& 0 & 0  & \cdots & 3  &  0\\
 1& 0 & 0  & \cdots & 0  &  1\\
  \end{array}\right).$$

Further,  we can write $B_1(H_s^k)$  as follows:
  $$B_1(H_s^k) = \left(\begin{array}{ccc}
 (n-2)J_1+I_1 & 2J_{1,s} & kJ_1\\
 J_{s,1} & 3I_s& 0 \\
 J_1 & 0 &  I_1\\
  \end{array}\right).$$
  Let $B_2(H_s^k)$ be  the corresponding equitable quotient matrix of $B_1(H_s^k)$. Then   by  Lemma \ref{lem29},  we have
$\sigma(B_1(H_s^k))=\sigma(B_2(H_s^k))\cup\{3^{[s-1]}\},$   where
  $$B_2(H_s^k) = \left(\begin{array}{ccc}
 n-1 &2s & k\\
 1 & 3& 0 \\
 1 & 0 &  1\\
  \end{array}\right).$$
Thus we have
\begin{equation}\label{eq25}
\sigma(Q(H_s^k))=\sigma(B_2(H_s^k))\cup\{1^{[\frac{n+k-3}{2}]},3^{[\frac{n-k-3}{2}]}\},
\end{equation}
and by direct computing, we know the signless Laplacian characteristic polynomial of $B_2(H_s^k)$ is as follows:
\begin{equation}\label{eq26}
\det(xI_n-B_2(H_s^k))=x^3-(n+3)x^2+3nx-2n+2k+2.\end{equation}
Therefore, we complete the proof of (\ref{eq23}) by (\ref{eq25}) and (\ref{eq26}).
\end{Proof}

\begin{lem}\label{lem211}
Let  $n, k$ be positive integers and $s$ be nonnegative integer with $2s+k+2=n$,
$L_s^k$ be a graph on $n$ vertices in  Fig. 2, and $f(x)=x^5-(n+5)x^4+(6n+4)x^3-(12n-2k-10)x^2+(9n-6k-12)x-2n+2k+4$. Then
\begin{equation}\label{eq27}
\psi(L_s^k)=(x-1)^{\frac{n+k-6}{2}}(x-3)^{\frac{n-k-4}{2}}f(x).
\end{equation}
\end{lem}

\begin{Proof}
Clearly, we have %We note that $2s+t+2=n$ and
%As the signless Laplacian matrix of $L_s^t$  as follows,
$$Q(L_s^k) = \left(\begin{array}{ccccccccccc}
 n-2 & 1 & 1 & \cdots & 1 & 1 & 1 & 0 & 1 & \cdots & 1\\
 1 & 2 & 1  & \cdots  & 0 & 0 & 0 & 0 & 0 & \cdots & 0\\
 1 & 1 & 2 & \cdots & 0 & 0 & 0 & 0 & 0 & \cdots & 0\\
 \vdots & \vdots & \vdots & \ddots  & \vdots &  \vdots & \vdots & \vdots & \vdots & \ddots & \vdots\\
 1 & 0 & 0 &\cdots  &2 & 1 & 0 & 0 & 0 & \cdots & 0\\
 1 & 0 & 0 & \cdots & 1 & 2 & 0 & 0 & 0 & \cdots & 0\\
 1 & 0 & 0 &\cdots & 0 &0 & 2 & 1 & 0 & \cdots & 0\\
 0 & 0 & 0 &\cdots& 0 & 0 & 1 & 1 & 0 & \cdots & 0\\
 1 & 0 & 0 &\cdots & 0 & 0 & 0 & 0 & 1 & \cdots & 0\\
 \vdots & \vdots & \vdots & \ddots  & \vdots &  \vdots & \vdots & \vdots & \vdots & \ddots & \vdots\\
 1 & 0 & 0 &\cdots  & 0 & 0 & 0 & 0 & 0 & \cdots & 1\\
  \end{array}\right)$$
and it can be written as follows:
  \begin{equation}\label{eq28}
  Q(L_s^k) = \left(\begin{array}{ccccccc}
 (n-3)J_1+I_1 & J_{1,2} & \cdots & J_{1,2}  & J_1& 0 &J_{1,k-1}\\
 J_{2,1} & J_2+I_2& \cdots  & 0 & 0 & 0  & 0\\
 \vdots & \vdots & \ddots  & \vdots &  \vdots & \vdots & \vdots\\
 J_{2,1} & 0 & \cdots & J_2+I_2 & 0 & 0& 0\\
J_1 &  0  &  \cdots &  0 & J_1+I_1& J_1& 0\\
 0& 0 &  \cdots & 0  &  J_1& I_1& 0\\
 J_{k-1,1}&0&\cdots & 0  &0 & 0  &  I_{k-1}\\
  \end{array}\right).\end{equation}

 Let $B_1(L_s^k)$ be  the corresponding equitable quotient matrix of $Q(L_s^k)$, then by (\ref{eq28}) and Lemma \ref{lem29},
 we have $\sigma(Q(L_s^k))=\sigma(B_1(L_s^k))\cup\{1^{[\frac{n+k-6}{2}]}\}$, where
 $$B_1(L_s^k) = \left(\begin{array}{ccccccc}
 n-2 & 2 &  \cdots & 2 &1 & 0 &k-1\\
 1 & 3 & \cdots & 0&0  &0  & 0\\
 \vdots &  \vdots  & \ddots &  \vdots & \vdots&  \vdots&  \vdots\\
 1& 0 & \cdots & 3  &  0 &  0 &  0\\
  1& 0 & \cdots & 0  &  2 &  1 &  0\\
  0& 0 & \cdots & 0  &  1 &  1&  0\\
   1& 0 & \cdots & 0  &  0 &  0 &  1\\
  \end{array}\right).$$

Further, we can write   $B_1(L_s^k)$ as follows:
  $$B_1(L_s^k) = \left(\begin{array}{ccccc}
 (n-3)J_1+I_1 & 2J_{1,s} & J_1 & 0 & (k-1)J_1\\
 J_{s,1} & 3I_{s}& 0  & 0 & 0\\
J_1 & 0& J_1+I_1  & J_1 & 0\\
0 & 0& J_1  &I_1 & 0\\
J_1 & 0& 0  & 0 &I_1\\
  \end{array}\right).$$

Let $B_2(L_s^k)$ be  the corresponding equitable quotient matrix of $B_1(L_s^k)$, then by Lemma \ref{lem29},
we have $\sigma(B_1(L_s^k))=\sigma(B_2(L_s^k))\cup\{3^{[s-1]}\},$  where
  $$B_2(L_s^k) = \left(\begin{array}{ccccc}
n-2 &2s& 1  & 0 &k-1\\
1 & 3& 0  & 0 & 0\\
1 & 0& 2  & 1 & 0\\
0 & 0& 1  & 1 & 0\\
1 & 0& 0 & 0 & 1\\
  \end{array}\right).$$

Thus
 \begin{equation}\label{eq29}
 \sigma(Q(L_s^k))=\sigma(B_2(L_s^k))\cup\{1^{[\frac{n+k-6}{2}]},3^{[\frac{n-k-4}{2}]}\},\end{equation}
and by direct computing, we know the signless Laplacian characteristic polynomial  of $B_2(L_s^k)$,
$\det(xI_n-B_2(L_s^k))=f(x).$

Combining the above arguments, (\ref{eq27}) holds.
\end{Proof}

\begin{prop}{\rm(\cite{2012LAA}, (3.1), (3.2))}\label{prop212}
 Let $h(x)=x^5-(k+9)x^4-(n-7k-32)x^3+(4n-14k-54)x^2-(4n-7k-40)x+n-k-8$.
  Then the following equations $(\ref{eq210})$ and $(\ref{eq211})$ are the equations  $(3.1)$ and $(3.2)$ in {\cite{2012LAA}}, respectively.
 \begin{equation}\label{eq210}
%{\rm(\cite{2012LAA}, (3.1))}
 \psi(H_s^k)=(x-1)^{\frac{n+k-3}{2}}(x-3)^{\frac{n-k-3}{2}}[x^3-(k+6)x^2-(n-4k-12)x+n-k-7].
\end{equation}
\begin{equation}\label{eq211}
 \psi(L_s^k)=(x-1)^{\frac{n+k-6}{2}}(x-3)^{\frac{n-k-4}{2}}h(x).
\end{equation}
\end{prop}

Comparing  (\ref{eq23}) with (\ref{eq210}),  and comparing (\ref{eq27}) with (\ref{eq211}),
we can find that they are different equations. The reason is that in the proof of Eq.(3.1) in \cite{2012LAA},
$A_3=B_{n-k-2}$ should be $A_3=B_{n-k-2}-2(x^2-4x+3)^{\frac{n-k-3}{2}}.$
Thus Theorems 3.1-3.3 and  Conjecture 3.4 in \cite{2012LAA} may be revised as follows.

%\emph For example, the  signless Laplacian matrix of $H_2^1$ and $L_2^1$ are as follows, respectively.
% $$Q(H_2^1) = \left(\begin{array}{cccccc}
%5 & 1& 1 & 1 & 1 & 1\\
%1& 2 & 1 & 0 & 0 & 0\\
%1& 1 & 2 & 0 & 0 & 0\\
%1 & 0 & 0 & 2& 1 & 0\\
%1 & 0 & 0 & 1& 2 & 0\\
%1 & 0 & 0 & 0& 0 & 1\\
%  \end{array}\right), \qquad Q(L_2^1) = \left(\begin{array}{cccccccc}
% 6 & 1 & 1 & 1 & 1 & 1 & 0 & 1\\
% 1 & 2 & 1 & 0 & 0 & 0 & 0 & 0\\
% 1 & 1 & 2 & 0 & 0 & 0 & 0 & 0\\
% 1 & 0 & 0 & 2 & 1 & 0 & 0 & 0\\
% 1 & 0 & 0 & 1 & 2 & 0 & 0 & 0\\
% 1 & 0 & 0 & 0 & 0 & 2 & 1 & 0\\
% 0 & 0 & 0 & 0 & 0 & 1 & 1 & 0\\
% 1 & 0 & 0 & 0 & 0 & 0 & 0 & 1\\
 % \end{array}\right).$$
%\emph By directed computing signless Laplacian characteristic polynomials  \emph{of} $H_2^1$ and $L_2^1$ are
%$$det(xI-Q(H_2^1))=(x-3)(x-1)^2(x-2)(x^2-7x+4),$$
%$$det(xI-Q(L_2^1))=(x-3)(x-1)^2(x-2)(x^4-11x^3+30x^2-22x+4).$$

\begin{prop}{\rm(\cite{2012LAA}, Theorem 3.1)}
Let $n$ be positive integer, $k, s$ be nonnegative integers with $k<n$ and $s=\lfloor\frac{n-k-1}{2}\rfloor,$
$G$ be a  cactus with $n$ vertices and  $k$ pendant vertices.

{\rm (i) }  If $n-k \equiv 1 \pmod{2}$, then $q_1(G) \leq q_1(H_s^k)$,  with equality if and only if $G \cong H_s^k,$ where $H_s^k$ is depicted in Fig. 1 and $q_1(H_s^k)$ is the largest root of the equation $x^3-(n+3)x^2+3nx-2n+2k+2=0.$

{\rm (ii) }   If $n-k \equiv 0 \pmod{2}$, then $q_1(G) \leq q_1(L_s^{k})$,  with equality if and only if $G \cong L_s^{k},$ where $L_s^{k}$ is depicted in Fig. 2 and $q_1(L_s^{k})$ is the largest root of the equation $x^5-(n+5)x^4+(6n+4)x^3-(12n-2k-10)x^2+(9n-6k-12)x-2n+2k+4=0$.
\end{prop}

\begin{prop}{\rm(\cite{2012LAA}, Theorem 3.2)}
Let $n$ be positive integer and $s$ be nonnegative integer  with $s=\lfloor\frac{n-1}{2}\rfloor,$
$G$ be a  cactus on $n$ vertices. Then

{\rm (i) } $q_1(G)\leq\frac{(n+2)+\sqrt{n^2-4n+12}}{2}$ for odd $n,$ and the equality holds if and only if $G\cong H_s^0.$

{\rm (ii) } $q_1(G)\leq\frac{(n+1)+\sqrt{n^2-2n+9}}{2}$ for even $n,$ and the equality holds if and only if $G\cong H_s^1.$
\end{prop}

\begin{prop}{\rm(\cite{2012LAA}, Theorem 3.3)}\label{prop215}
Let $n$ be even,   $G$ be a  cactus with $n$ vertices and a perfect matching. Then $q_1(G) \leq \frac{(n+1)+\sqrt{n^2-2n+9}}{2},$ and the equality holds if and only if $G \cong H_{m-1}^1,$ where $m=\frac{n}{2}$.
\end{prop}

\begin{con}{\rm(\cite{2012LAA}, Conjecture 3.4)}\label{con216}
Let $G$ be a cactus on $n$ vertices with  matching number $m.$

{\rm (i) } If $n=2m+1,$ then $q_1(G)\leq\frac{(n+2)+\sqrt{n^2-4n+12}}{2},$ and the equality holds if and only if $G \cong H_m^0.$

{ \rm (ii) } If $n\geq2m+2,$ then $q_1(G) \leq q_1(H_{m-1}^{n-2m+1}),$  and the equality holds if and only if $G \cong H_{m-1}^{n-2m+1},$ where $q_1(H_{m-1}^{n-2m+1})$ is the largest root of the equation $x^3-(n+3)x^2+3nx-4m+4=0.$
\end{con}

\section{Main result}

%\vskip.15cm
\hskip.6cm  Let  $G$ be  a cactus  with $n$ vertices and matching number $m.$ Clearly, $m\leq \frac{n}{2}$.
If $n=2m$, Proposition \ref{prop215} characterized the unique graph with the maximum signless Laplacian spectral radius
among all cacti in $\ell_n^m$.
%gave the sharp upper bound of the signless Lpalacian spectral radius of cacti   with $n$ vertices and matching number $m.$
Now we will  determine the case of $n\geq 2m+1$ and thus show Conjecture \ref{con216} is true.
The technique used in the proof is motivated by \cite{2014,2012LAA} et al.
\begin{them}\label{thm31}
Let $G$ be a graph in $\ell_n^m.$

{\rm (i) } If $n=2m+1,$ then $q(G)\leq\frac{(n+2)+\sqrt{n^2-4n+12}}{2}$  with equality if and only if $G \cong H_m^0.$

{\rm (ii) } If $n\geq2m+2,$ then $q(G) \leq q(H_{m-1}^{n-2m+1})$ with equality if and only if $G \cong H_{m-1}^{n-2m+1},$ where $q(H_{m-1}^{n-2m+1})$ is the largest root of the equation of $x^3 - (n+3)x^2 + 3nx - 4m +4 =0.$
\end{them}

\begin{Proof}
Let $G=(V,E) \in \ell_n^m$ such that its signless Laplacian spectral radius is the largest in $\ell_n^m$,
where $V=\{v_1,v_2,\ldots,v_n\}$. Let  $\textbf{x}=(x_{v_1},x_{v_2},\ldots,x_{v_n})^T$ be the Perron vector of $Q(G)$,
 where $x_{v_i}$ corresponds to the vertex $v_i, i = 1,2,\ldots,n.$

Let $M$ be a maximum matching of $G.$ Then $|M|=m,$ and there are three cases for a non-pendant edge $e=uv$ in $G:$
(1) $e=uv$ is an $M$-saturated edge; (2) $e=uv$ has exactly one $M$-saturated vertex;
(3) $e=uv$ is  an $M$-unsaturated edge but both $u$ and $v$ are $M$-saturated vertices.

We first prove that the graph $G$ is a bundle. In order to do so we will prove the following  four claims.

\vskip.15cm
\textbf{Claim 1.}  If $T$ is an induced subtree of $G$, and $T$ is  attached to a vertex $v$, where $v$ is on some cycle of $G,$
then $T$ is a star whose center is $v$.
\begin{proof}
Suppose that there exists a  tree $T$ attached to $v_{1},$ but $T$ is not a star, where $v_1$ is on some cycle of $G$.
Then there must exist a path $P=v_{1}u_1u_2$ of length 2 in $T.$ Take $V_1=N_G(u_1)\backslash \{v_{1}\}$, $V_2=N_G(v_{1})\backslash \{u_1\}$ and
$$G_1=\left\{\begin{array}{cc}
 G-\{u_1x\hskip.1cm |\hskip.1cm x\in V_1\}+\{v_{1}x\hskip.1cm |\hskip.1cm x\in V_1\}, & \mbox { if } x_{v_{1}}\geq x_{u_1}; \\
 G-\{v_{1}x\hskip.1cm |\hskip.1cm x\in V_2\}+\{u_1x\hskip.1cm |\hskip.1cm x\in V_2\}, & \mbox { if } x_{v_{1}}< x_{u_1}.
                \end{array}\right.$$

In $G_1$, $T$ is transformed to a new tree, denoted by $T_1$.
It is obvious that the  number of pendant vertices of $T_1$ is more than the number of pendant vertices of $T$.
If there still exists a path $Q$ of length 2 in $T_1$ of $G_1$, repeatedly the similar operation as above,
till there are no  path with length 2 in $T_1$. At this moment, $T_1$ is transformed to a star, denoted by $T_2$,
and the new graph is denoted by $G_2.$  By Lemma \ref{lem21}, we have $q(G)<q(G_1)<q(G_2)$.

 Assume that $|V(T)|=n_1$, so there are at most  $\lfloor\frac{n_1}{2}\rfloor $ edges of $T$ in the maximum matching of $G$.
 Therefore,  we can add a new edge   by connecting two pendant vertices of $T_2$, and add another new edge by connecting another two pendant vertices of $T_2$, till the new graph, say $G^*$,  satisfying $G^*\in \ell_n^m.$
By Corollary \ref{cor26}, we have $q(G_2)<q(G^*)$, and thus $q(G)<q(G^*),$  it is a contradiction by the assumption of $G$.
\end{proof}

\textbf{Claim 2.} All cycles of $G$ are with length 3.
\begin{proof}
We assume that there exists a cycle, say $C^0={v_1 v_2 \ldots v_l v_1}$ in $G,$ such that its length  $l\geq 4$.
 As $G$ is a cactus, we know $v_i v_j \not\in E_G$ for any $|i-j|\geq 2$. We will complete the proof of Claim 2 by the following three cases.

\textbf{Case 1:} There exists an edge $uv$ on $C^0$  such that $uv \in M .$

It is obvious that   $N(u) \cap N(v) = \emptyset$ by the definition of a cactus.
Let $G_1$ be the graph obtained from $G$ by deleting  edge $uv$, identifying vertex $u$ and vertex $v$,
adding a new vertex $w$ and a new pendant edge $uw$ to $u$.
By Lemma \ref{lem22}, $q(G)<q(G_1).$
Take $M_1= M-uv+uw$, then $M_1$ is a matching of $G_1.$ It is easy to see that $M_1$ is a maximum matching of $G_1,$
thus $G_1 \in \ell_n^m,$  it is a contradiction by the  assumption of $G$.

\textbf{Case 2:} There exists an edge $uv \in C^0$  having exactly one $M$-saturated vertex.

Without loss of generality, we assume that $u$ is  an $M$-saturated vertex and $v$ is an $M$-unsaturated vertex.
Similar to the proof of Case 1, we obtain $G_1$ from $G$ and
take $M_1= M$. Then we can show $M_1$ is  a maximum matching of $G_1$, thus $G_1 \in \ell_n^m,$ and  $q(G)<q(G_1),$  it is a contradiction.

\textbf{Case 3:}  All edges on $C^0$ are  $M$-unsaturated edges but all vertices of $C^0$ are $M$-saturated.

Then for any $i\in\{1,2,\ldots, l\}$, $v_i$ is adjacent to a vertex  $w_i$
such that  $v_iw_i \in M$ and %each $v_iw_i$ is a pendant edge by Claim 1.
$w_i\neq v_{i-1}, w_i \neq v_{i+1}$ (if $i=l,v_{i+1}=v_1$; if $i=1, v_{i-1}=v_l$).
%Furthermore,  for each $i,j\in \{1,2,\ldots, l\}$,  there exists exactly one path $w_iv_i\ldots v_jw_j$ from
%$w_i$ to $w_j$ by the definition of a cactus.

\textbf{Subcase 3.1:}  $x_{v_1} \geq x_{v_2}.$

Let $G_1=G - v_2 v_3 + v_1 v_3.$  Then  we have $q(G)<q(G_1)$  by Lemma \ref{lem21}.
Take $M_1= M$, then $M_1$ is a matching of $G_1.$
Now we show that  $M_1$ is a maximum matching of $G_1$.

If $M_2$ is a maximum matching of $G_1$ and $|M_2|>m$, then $v_1v_3\in M_2$
(otherwise, $M_2\subset E(G_1-v_1v_3)\subset E(G)$ is a matching of $G$, it is a contradiction).
If $w_1, w_3$ are $M_2$-saturated,  then we take $M_3=M_2-v_1v_3$. Clearly, $M_3$ is a matching of $G$ with $|M_3|=|M_2|-1\geq m$.
If $|M_3|>m$, it implies a contradiction by $G \in \ell_n^m$. Thus $|M_3|=m$ and $M_3$ is also a maximum matching of $G$.
We note that $v_1$ is $M_3$-unsaturated, then $v_2$ must be  $M_3$-saturated. Similar to the proof of Case 2, we can obtain a contradiction.
Therefore $w_1$ is $M_3$-unsaturated or $w_3$ is $M_3$-unsaturated. Without loss of generality, we assume that $w_1$ is $M_3$-unsaturated.
Then $M_4=M_2-v_1v_3+v_1w_1$ is a matching of $G$ with $|M_4|=|M_3|>m$, it is   a contradiction by $G \in \ell_n^m$.
Combining the above arguments, we have $G_1 \in \ell_n^m,$ it is a contradiction by the assumption of $G$.

\textbf{Subcase 3.2:}  $x_{v_1} < x_{v_2}.$

Let $G_1=G - v_1 v_l + v_2 v_l.$ Similar to the proof of Subcase 3.1, we can show  a contradiction, so we omit it.
\end{proof}

\textbf{Claim 3. } Any two cycles of $G$ have exactly one common vertex.
\begin{proof}
On the contrary, by  any two cycles of $G$ have no common edges,
we can  assume  that there are two disjoint cycles $C^1$ and $C^2$  in $G$ connecting by
 the shortest path $P=v_1 v_2 \ldots v_s$ of length $s-1(s \geq 2)$  with $V_P \cap V_{C^1}=\{v_1\}$ and $V_P \cap V_{C^2}=\{v_s\}.$
  In view of Claim 2, let $C^1=v_1 w_1 w_2 v_1,C^2=v_s u_1 u_2 v_s.$
  Now we will show that there exists an $m$-matching, say $M_1,$ such that $\{w_1 v_1, w_2 v_1, u_1 v_s, u_2 v_s\} \cap M_1 =\emptyset.$
  Noting that the proof of  $\{u_1v_s,u_2v_s\} \cap M_1=\emptyset$ is similar to the proof of  $\{w_1 v_1,w_2 v_1\} \cap M_1 =\emptyset,$
    now we only show   $\{w_1 v_1,w_2 v_1\} \cap M_1 =\emptyset.$

In fact, if ${w_1 v_1,w_2 v_1}\notin M,$ we take $M_1=M,$ then our result holds. Otherwise, without loss of generality, we assume that $w_2 v_1\in M.$

{\bf Case 1: } $w_1$ is not saturated by $M.$

Let $M_1=M - w_2 v_1 + w_1 w_2.$ It is easy to see that $M_1$ is an $m$-matching of $G$ and $\{w_1 v_1,w_2 v_1\} \cap M_1 =\emptyset,$
 then our result holds in this case.

{\bf Case 2: } $w_1$ is saturated by $w_1 t_1$ in $M.$

In this case, $d_G(w_1)\geq 3.$ Let
$V_1=N_G(v_1) \backslash \{w_1 ,w_2\}$, $V_2=N_G(w_1) \backslash \{v_1 ,w_2\}$ and
$$G_1=\left\{\begin{array}{cc}
 G-\{v_1x\hskip.1cm |\hskip.1cm x \in V_1\} + \{w_1x\hskip.1cm |\hskip.1cm x \in V_1\}, & \mbox { if }  x_{w_1} \geq x_{v_1}; \\
 G-\{w_1x\hskip.1cm |\hskip.1cm  x \in V_2\} + \{v_1x\hskip.1cm |\hskip.1cm x \in V_2\}, &  \mbox { if } x_{w_1} < x_{v_1}.
            \end{array}\right.$$
 By Lemma \ref{lem21}, $q(G)<q(G_1).$ Take
 $$M_2=\left\{\begin{array}{cc}
 M, & \mbox { if }  x_{w_1} \geq x_{v_1}; \\
 M- w_1 t_1 - w_2 v_1 + w_1 w_2 + v_1 t_1, &  \mbox { if } x_{w_1} < x_{v_1}.
            \end{array}\right.$$
Then  $M_2$ is an $m$-matching of $G_1.$

Now we  show  $G_1 \in \ell_n^m.$ Otherwise, let $M_3$ be a maximum matching of $G_1$ with $|M_3|>m$,
now we will obtain a  contradiction by the following two subcases.

{\bf Subcase 2.1: } $x_{w_1} \geq x_{v_1}$.

Then $w_1t\in M_3$ where $t\in N_G(v_1)\backslash \{w_1,w_2\}$.
Otherwise, $M_3\subset E(G_1)-\{w_1t\hskip.05cm |\hskip.05cm t\in N_G(v_1)\backslash\{w_1,w_2\}\}\subset E(G)$ is a matching of $G$, it is a contradiction by $G\in \ell_n^m$ and $|M_3|>m$.

If $v_1$ is $M_3$-saturated, then $v_1w_2\in M_3$ by $d_{G_1}(v_1)=2$. Thus we take
 $$M_4=\left\{\begin{array}{cc}
                          M_3-w_1t-v_1w_2+v_1t+w_1w_2, &  v_1 \mbox{ is } M_3\mbox{-saturated};\\ %v_1w_2\in M_3; \\
                         M_3-w_1t+v_1t, & v_1 \mbox{ is } M_3\mbox{-unsaturated}.
                        \end{array}\right.$$
                               \noindent Clearly, $M_4$ is  an $|M_3|$-matching of $G$, it is a contradiction by $G\in \ell_n^m$ and $|M_3|>m$.

{\bf Subcase 2.2: } $x_{w_1} < x_{v_1}$.

Then $v_1t\in M_3$ where $t\in N_G(w_1)\backslash\{v_1,w_2\}$.
Otherwise, $M_3\subset E(G_1)-\{v_1t\hskip.05cm |\hskip.05cm t\in N_G(w_1)\backslash\{v_1,w_2\}\}\subset E(G)$ is a matching of $G$,
it is a contradiction by $G\in \ell_n^m$ and $|M_3|>m$.

If $w_1$ is $M_3$-saturated, then $w_1w_2\in M_3$ by $d_{G_1}(w_1)=2$. Take
 $$M_4=\left\{\begin{array}{cc}
                          M_3-v_1t-w_1w_2+w_1t+v_1w_2, &  w_1 \mbox{ is } M_3\mbox{-saturated};\\ %w_1w_2\in M_3; \\
                         M_3-v_1t+w_1t, & w_1 \mbox{ is } M_3\mbox{-unsaturated}.
                        \end{array}\right.$$
                               \noindent Clearly, $M_4$ is  an $|M_3|$-matching of $G$, it is a contradiction by $G\in \ell_n^m$ and $|M_3|>m$.

Combining the above two subcases, we have $G_1 \in \ell_n^m.$
Then it implies a contradiction by the assumption of $G$, $q(G)<q(G_1)$ and $G_1 \in \ell_n^m.$

Thus, there exists an $m$-matching $M_1$ of $G$ such that $\{w_1 v_1,w_2 v_1,u_1 v_s,u_2 v_s\} \cap M_1 =\emptyset.$
Now we take
$$G_2=\left\{\begin{array}{cc}
             G-v_1 w_1-v_1 w_2+v_s w_1+v_s w_2, & \mbox {if }  x_{v_s} \geq x_{v_1}; \\
              G-v_s u_1-v_s u_2+v_1 u_1+v_1 u_2, & \mbox {if }  x_{v_s} < x_{v_1}.
            \end{array}\right.$$

 By Lemma \ref{lem21}, $q(G)<q(G_2).$ Let $M_5=M_1$, it is obvious that $M_5$ is an $m$-matching of $G_2$.

 Now we show  $G_2 \in \ell_n^m.$  % which it implies  a contradiction by the assumption of $G$.
Otherwise, let   $M_6$ be a maximum matching of $G_2$ with $|M_6|>m$.
 We only consider the case  $x_{v_s} \geq x_{v_1}$ because the proof of the case $x_{v_s} < x_{v_1}$ is similar to the case $x_{v_s} \geq x_{v_1}$.

 Then $v_sw_1\in M_6$ or $v_sw_2\in M_6$. Otherwise $M_6\subset E(G_2-v_s w_1-v_s w_2)\subset E(G)$ is a matching of $G$,
 it is a contradiction by $G\in \ell_n^m$ and $|M_6|>m$.
 Without loss of generality, we assume that $v_sw_1\in M_6$.
 Then we have the following results.

  (i) There exists some  $t\in N_G(v_1)\backslash\{w_1, w_2\}$ such that $v_1t\in M_6$.
  Otherwise, $M_6-v_sw_1+v_1w_1\subset E(G)$ is an $|M_6|$-matching of $G$, it is a contradiction by $G\in \ell_n^m$ and $|M_6|>m$.

 (ii) There exists some  $r\in N_G(w_2)\backslash\{v_1, w_1\}$ such that $w_2r\in M_6$.
  Otherwise, $M_6-v_sw_1+w_1w_2\subset E(G)$ is an $|M_6|$-matching of $G$, it is a contradiction by $G\in \ell_n^m$ and $|M_6|>m$.

 (iii) There exist no pendant vertices $z$ such that $w_1z\in E(G)$.
 Otherwise, $M_6-v_sw_1+w_1z\subset E(G)$ is an $|M_6|$-matching of $G$, it is a contradiction by $G\in \ell_n^m$ and $|M_6|>m$.

By (ii), we have $d_{G}(w_2)\geq 3$. Then we let $V_3=N_G(w_2)\backslash\{v_1,w_1\}$ and

 $$G_3=\left\{\begin{array}{cc}
 G-\{v_1x\hskip.1cm |\hskip.1cm x \in V_1\} + \{w_2x\hskip.1cm |\hskip.1cm x \in V_1\}, & \mbox { if }  x_{w_2} \geq x_{v_1}; \\
 G-\{w_2x\hskip.1cm |\hskip.1cm  x \in V_3\} + \{v_1x\hskip.1cm |\hskip.1cm x \in V_3\}, &  \mbox { if } x_{w_2} < x_{v_1}.
            \end{array}\right.$$
            By Lemma \ref{lem21}, $q(G)<q(G_3).$
Similar to the proof of Case 2, we can show  $G_3\in \ell_n^m$, it is  a contradiction by the assumption of $G$.

By the above arguments, we know   $G_2 \in \ell_n^m$.
But it implies a contradiction by the assumption of $G$, $q(G)<q(G_2)$ and $G_2 \in \ell_n^m.$
Combining the above arguments, we  completes the proof of Claim 3.
\end{proof}

\textbf{Claim 4.} Any three cycles contained in $G$ have exactly one common vertex.
\begin{proof}
We assume that there exists three cycles say $C^1,C^2$ and $C^3$ in $G$ such that they have no common vertex.
By Claim 2, let $V_{C^1}\cap V_{C^2}=\{u\},V_{C^1}\cap V_{C^3}=\{v\}$ and $ V_{C^2}\cap V_{C^3}=\{w\}.$
Then $u,v$ and $w$ at the same cycle $C^4.$ Thus $C^1\cap C^4=\{u,v\},$ a contradiction to the definition of a cactus.
Then any three cycles have  common vertices, and thus any three cycles  have exactly one common vertex by Claim 3.
\end{proof}

By Claims 3 and 4, we know that all of the cycles contained in $G$ have exactly one common vertex, say $v_0,$ $ i.e. $ $G $ is a $bundle.$
Now we know that the graph in $\ell_n^m$ having the largest signless Laplacian spectral radius  is a $bundle$ with some pendant trees attached,
 and these trees are stars by Claim 1. Now we will show these stars are attached to the vertex $v_0$ of $G$.

\vskip.15cm
\textbf{Claim 5.} If $G$ contains pendant edges, then all pendant edges are attached to the common vertex $v_0$.

\begin{proof}
On the contrary,  by Claim 2 we can assume that there exists a cycle, say $v_0v_1v_2v_0,$
in $G$ such that $v_1$ is adjacent to $i$  pendant vertices and $v_2$ is adjacent to $j$  pendant vertices,
where $i,j\geq 0$ and $i+j\geq 1.$ We will complete the proof  by the following four cases.

\textbf{Case 1:} $i\geq 1,j=0$ or $i=0,j\geq 1.$

Without loss of generality, let $i\geq 1,j=0.$ The proof of the case $i=0,j\geq 1$ is similar, we omit it.
Let $V_1=N_G(v_1)\backslash\{v_0,v_2\}$, $V_2=N_G(v_0)\backslash\{v_1,v_2\}$.

\textbf{Subcase 1.1:} There exists a pendant edge, say $v_1 v'_1,$ attached to $v_1$ being in $M.$

 Take
\begin{equation}\label{eq31}
G_1=\left\{\begin{array}{cc}
              G-\{v_1x\hskip.1cm | \hskip.1cm x \in V_1\}
              +\{v_0x \hskip.1cm | \hskip.1cm x\in V_1\},  & \mbox { if } x_{v_0}\geq x_{v_1};\\
              G-\{v_0x\hskip.1cm | \hskip.1cm x \in V_2\}
              +\{v_1 x \hskip.1cm | \hskip.1cm x\in V_2\}, & \mbox { if } x_{v_0} < x_{v_1}.
            \end{array}\right.\end{equation}
  By Lemma \ref{lem21}, $q(G)<q(G_1).$

  In the case $x_{v_0}\geq x_{v_1}$, we set
  $$M_1=\left\{\begin{array}{cc}
                                                           M-v_1 v'_1-v_0 v_2+v_1 v_2+v_0 v'_1, & \mbox { if } v_0v_2\in M; \\
                                                          M - v_1 v'_1 + v_1 v_2, &  \mbox { if } v_0v_2\not\in M.
                                                        \end{array}\right.$$
  Then $M_1$ is an $m$-matching of $G_1.$ Now we show $M_1$ is a maximum matching of $G_1$.
  If $M_2$ is a maximum matching of $G_1$ and $|M_2|>m$, then there exists $u\in V_1$ such that $v_0u\in M_2$.
  If not, $M_2\subset E(G_1-\{v_0x \hskip.1cm | \hskip.1cm x \in V_1\})\subset E(G)$ is a  matching of $G$, it is a contradiction.
  What's more,  we have $v_1v_2\in M_2$ by $d_{G_1}(v_1)=d_{G_1}(v_2)=2$.
  Then $M_3=M_2-v_1v_2-v_0u+v_1u+v_0v_2$ is a   matching of $G$ with $|M_3|=|M_2|>m$, it is a contradiction by $G \in \ell_n^m$.
  Therefore $G_1 \in \ell_n^m,$  and thus it is a contradiction by the assumption of $G$.

  In the case $x_{v_0}<x_{v_1}$, we set
  $$M_1=\left\{\begin{array}{cc}
                                                           M, & \mbox { if } v_0v_2\in M; \\
                                                          M - v_0u + v_0 v_2, &  \mbox { if } v_0v_2\not\in M,
                                                        \end{array}\right.$$
where $u\in V_2$ and  $uv_0\in M$.   Then $M_1$ is  an  $m$-matching of $G_1,$
  and it is easy  to check that $G_1 \in \ell_n^m,$  thus it is a contradiction  by the assumption of $G$.

\textbf{Subcase 1.2:} Each of the pendant edges attached to $v_1$ is not in $M.$

Then $v_1 v_0\in M$ or $v_1 v_2\in M.$
If $v_1 v_0\in M,$  we set $M_1=M-v_1v_0+v_0v_2+v_1u$, where $u\in V_1$.
Then $M_1$ is an $(m+1)$-matching of $G,$  it implies a contradiction by the assumption of $G.$

Thus $v_1 v_2 \in M.$ Then $v_0$ must be $M$-saturated by $v_0w$ where $w\in V_2$. % and $w$ must be a pendant vertex.
Otherwise, if $v_0$ is not $M$-saturated, we take $M_1=M-v_1v_2+v_1u+v_0v_2$ where $u\in V_1$ is an $(m+1)$-matching of $G$,
it is a contradiction by the assumption of $G.$ Thus  $v_0$ must be $M$-saturated by $v_0w$ where $w\in V_2$.

Take $G_1$ as (\ref{eq31}), then  $q(G)<q(G_1)$ by Lemma \ref{lem21}. We set
$$M_1=\left\{\begin{array}{cc}
                                                           M, & \mbox { if } x_{v_0}\geq x_{v_1}; \\
                                                          M-v_1 v_2-v_0w+v_1w+v_0v_2, &  \mbox { if } x_{v_0}< x_{v_1}.
                                                        \end{array}\right.$$
  Then $M_1$ is an  $m$-matching of $G_1.$
 Similar to the proof of Subcase 1.1,  it is easy  to check that $G_1 \in \ell_n^m,$  thus it is a contradiction  by the assumption of $G$.

\textbf{Case 2:} $i\geq 2,j\geq 2.$

If $x_{v_1}\geq x_{v_2},$ it is easy to see that there exists a pendant edge $v_2v'_2,$ such that $v_2v'_2\notin M.$
Let $G_1=G-v_2v'_2+v_1v'_2,$ By Lemma \ref{lem21}, $q(G)<q(G_1).$ Obviously,  $G_1 \in \ell_n^m,$  a contradiction.

If $x_{v_1} < x_{v_2},$ it is easy to see that there exists a pendant edge $v_1v'_1,$ such that $v_1v'_1\notin M.$
Let $G_1=G-v_1v'_1+v_2v'_1.$ By Lemma \ref{lem21}, $q(G)<q(G_1).$ Obviously,  $G_1 \in \ell_n^m,$  a contradiction.

\textbf{Case 3:} $i\geq 2,j=1$ or $i=1,j\geq 2.$

Without loss of generality, we assume that $i\geq 2, j=1.$
Let $v_2v'_2$ be the only pendant edge attached to $v_2$.
It is easy to see that there exists a pendant edge $v_1v'_1$ such that $v_1v'_1\notin M$ by $i\geq 2$.

{\bf Subcase 3.1: }  $x_{v_1} < x_{v_2}.$

Let $G_1=G-v_1 v'_1+v_2v'_1.$ By Lemma \ref{lem21}, $q(G)<q(G_1).$ Obviously, $G_1 \in \ell_n^m,$ a contradiction.

{\bf Subcase 3.2: }   $x_{v_1} \geq x_{v_2}.$ Noting that $\{v_2'\}=N_G(v_2)\setminus \{v_0,v_1\},$  we let $G_1=G-v_2 v'_2+v_1v'_2.$ By Lemma \ref{lem21}, $q(G)<q(G_1).$

\textbf{Subcase 3.2.1: } $v_2 v'_2 \notin M.$

Then $M_1=M$ is  an $m$-matching of $G_1,$ and it is obvious to show  $G_1 \in \ell_n^m,$
which implies  a contradiction by the assumption of $G.$

\textbf{Subcase 3.2.2: } $v_2 v'_2 \in M.$

Then the matching number of $G_1$ is less than or equal to $m$. If not, there exists a matching $M_1$ of $G_1$ with $|M_1|>m$,
then $v_1v'_2\in M_1$. Otherwise, if $v_1v'_2\not\in M_1,$ $M_1\subset E(G_1-v_1v'_2)\subset E(G)$, and thus $M_1$ is a matching of $G$ with $|M_1|>m=|M|$, which implies  a contradiction by  $G\in \ell_n^m.$
Therefore, $M_2=M_1-v_1v'_2+v_1v'_1$ is a matching of $G$ with $|M_2|=|M_1|>m$, which implies  a contradiction by  $G\in \ell_n^m.$

If the matching number of $G_1$ is  equal to $m$, say, $G_1\in \ell_n^m,$ we obtain    a contradiction by the assumption of $G.$

 If the matching number of $G_1$ is  less than $m$, we set $G_2=G_1+v'_1v'_2,$ then $q(G)<q(G_1)<q(G_2)$ by Lemma \ref{lem21} and Corollary \ref{cor26}. On the other hand, we will show $G_2\in \ell_n^m,$  which implies  a contradiction by the assumption of $G.$

First, $M-v_2 v'_2+v'_1v'_2$ is  an $m$-matching of $G_2$. Second, if there exists a matching $M_3$ of $G_2$ with $|M_3|>m$,
then $v'_1v'_2\in M_3$. Otherwise, if $v'_1v'_2\not\in M_3,$ $M_3\subset E(G_2-v'_1v'_2)=E(G_1)$, and thus $M_3$ is a matching of $G_1$ with $|M_3|>m$, which implies  a contradiction because the matching number of $G_1$ is  less than $m$.
Further, $M_4=M_3-v'_1v'_2\subset E(G_2-v'_1v'_2)=E(G_1)$, and thus $M_4$ is a matching of $G_1$ with $|M_4|\geq m$, which implies  a contradiction because the matching number of $G_1$ is  less than $m$.
Combining the above arguments, $G_2\in \ell_n^m$ and we complete the proof of Case 3.

\textbf{Case 4:}   $i=j=1.$

Let $v_1v'_1$  be the only pendant edge attached to $v_1$,
and $v_2v'_2$  be the only pendant edge attached to $v_2$.
Then $v_1v'_1, v_2v'_2$ at least one in $M.$ Otherwise, $v_1v'_1\notin M $ and $v_2v'_2\notin M,$ then $v_1v_2\in M.$
Now we  take $M_1=M-v_1v_2+v_1v'_1+v_2v'_2$. Obviously, $M_1$ is an  $(m+1)$-matching of $G$,
it is  a contradiction by the assumption of $G.$

\textbf{Subcase 4.1:} $v_1 v'_1\in M,v_2 v'_2 \notin M $ or $v_1 v'_1 \notin M,v_2v'_2 \in M.$

 Without loss of generality, we assume $v_1 v'_1\in M$ and $v_2 v'_2 \notin M. $
Then $v_0v_2\in M.$ Otherwise, $M_1=M+v_2v'_2$   is an $(m+1)$-matching of $G,$
it is implies a contradiction by the assumption of $G.$

Take $G_1=\left\{\begin{array}{cc}
                   G-v_2 v'_2+v_1v'_2, & \mbox { if } x_{v_1} \geq x_{v_2};\\
                   G-v_1 v'_1+v_2v'_1, & \mbox { if } x_{v_1} < x_{v_2}.
                 \end{array}\right.$
 Then $q(G)<q(G_1) $ by Lemma \ref{lem21}.
  Further, we set $$M_1=\left\{\begin{array}{cc}
                  M, & \mbox { if } x_{v_1} \geq x_{v_2};\\
                  M-v_1v'_1-v_0v_2+v_0v_1+v_2v'_1, & \mbox { if } x_{v_1} < x_{v_2}.
                 \end{array}\right.$$ Then $M_1$ is  an  $m$-matching of $G_1,$  and it is obvious the check that $G_1 \in \ell_n^m,$
             which implies   a contradiction by the assumption of $G.$

\textbf{Subcase 4.2:} $v_1v'_1\in M$ and $v_2v'_2\in M.$

It is easy to see that $v_0v_1\not\in M$, $v_0v_2\not\in M$ and $v_1v_2\not\in M.$

\textbf{Subcase 4.2.1: }  $x_{v_1} \geq x_{v_2}.$

Let $G_1=G-v_2 v'_2+v_1v'_2,$ $G_2=G_1+v'_1v'_2.$ By Lemma \ref{lem21} and Corollary \ref{cor26}, we have $q(G)<q(G_1)<q(G_2).$
It is clear that $M_1=M-v_1v'_1-v_2v'_2 +v_1v_2+v'_1v'_2$  is an $m$-matching of $G_2.$
Now we show $G_2\in \ell_n^m.$

On the contrary, we assume that $M_2$ is a maximum matching of $G_2,$ and $|M_2|>|M|=m.$

 If $v'_1v'_2\in M_2,$ then $v_0v_1\in M_2,$ or $v_0v_2\in M_2,$ or $v_1v_2\in M_2.$
 Noting that $d_{G_2}(v_2)=2$ and $d_{G_2}(v_1)=4$, we set
 $$M_3=\left\{\begin{array}{cc}
                 M_2-v'_1v'_2-v_0v_1+v_1v'_1+v_2v'_2,  & \mbox { if } v_0v_1\in M_2; \\
                M_2-v'_1v'_2-v_0v_2+v_1v'_1+v_2v'_2,  & \mbox { if } v_0v_2\in M_2;\\
                M_2-v'_1v'_2-v_1v_2+v_1v'_1+v_2v'_2,  & \mbox { if } v_1v_2\in M_2.
               \end{array}\right.$$
               Clearly, $M_3$ is a matching of $G$  and $|M_3|= |M_2|>m,$  it is a contradiction by the assumption of $G.$
Thus $v'_1 v'_2\notin M_2.$ Therefore $v_1 v'_1\in M_2$ or $v_1 v'_2\in M_2.$

 If $v_1 v'_1\in M_2,$  we note that $M_2\subset E(G)$, then
 $M_2$ is a matching of $G.$ But $|M_2|>m,$ it is a contradiction by the assumption of $G.$

 If $v_1 v'_2\in M_2,$ then $v_0 v_2\in M_2.$
  Otherwise, we set $M_3=M_2-v_1v'_2+v_1v_2+v'_1v'_2$ is also a matching of $G_2,$
 but $|M_3|>|M_2|,$  it is a contradiction by he assumption of $M_2.$
Thus $M_4=M_2-v_0 v_2-v_1 v'_2+v_1 v'_1+v_2v'_2$ is a matching of $G,$
and $|M_4|>m,$ it is a contradiction by the assumption of $G.$

Combining the above arguments, we have  $G_2\in \ell_n^m,$  it is a contradiction by the assumption of $G.$

\textbf{Subcase 4.2.2: }   $x_{v_1} < x_{v_2}.$

Let $G_1=G-v_1v'_1+v_2v'_1,$ $G_2=G_1+v'_1v'_2.$ By Lemma \ref{lem21} and Corollary \ref{cor26}, we have $q(G)<q(G_1)<q(G_2).$
It is clear that $M_1=M-v_1v'_1-v_2v'_2 +v_1v_2+v'_1v'_2$  is  an  $m$-matching of $G_2,$
and similar to the proof of Subcase 4.2.1 we can  show $G_2\in \ell_n^m,$ it implies a contradiction by the assumption of $G.$
\end{proof}

Now we come back to complete the proof of Theorem 3.1.

By Claim 2, we conclude that the lengths of all cycles in $G$ are 3. By Claim 3 and Claim 4, all cycles have exactly one common vertex, say $v_0.$ By Claim 1 and Claim 5, all trees attached to a vertex of a cycle are stars (thus all trees are pendant edges) and these pendant edges attached to the same vertex $v_0.$ Therefore, if $n=2m+1,$ then $G\cong H_m^0$ or $G\cong H_{m-1}^2,$ where $H_m^0$ and $H_{m-1}^2$ are depicted in Fig. 1.
On the other hand,  by Corollary \ref{cor26}, we know that $q(H_{m-1}^2)<q( H_m^0),$  then $G\cong H_m^0$ when $n=2m+1.$
If $n\geq 2m+2,$ then $G\cong H_{m-1}^{n-2m+1},$ where $H_{m-1}^{n-2m+1}$ is depicted in Fig. 1.

By (\ref{eq23}), we have%$\psi(H_s^t)=(x-1)^{\frac{n+t-3}{2}}(x-3)^{\frac{n-t-3}{2}}[x^3-(n+3)x^2+3nx-2n+2t+2]$%Hence, we have
$$\psi(H_m^0)=(x-1)^{\frac{n-1}{2}}(x-3)^{\frac{n-3}{2}}[x^2-(n+2)x+2n-2]$$
and
$$\psi(H_{m-1}^{n-2m+1})=(x-1)^{n-m-1}(x-3)^{m-2}[x^3-(n+3)x^2+3nx-4m+4].$$

If $n\geq 3$ and $s\geq 1$, we note  that $C_3$ is a  subgraph of $H_s^t$, then  $q( H_s^t)\geq q(C_3)=4$  by Corollary \ref{cor26},
and thus $q( H_m^0)=\frac{(n+2)+\sqrt{n^2-4n+12}}{2}$ and $q( H_{m-1}^{n-2m+1})$ is the largest root of the equation $x^3-(n+3)x^2+3nx-4m+4=0.$
%This completes the proof of Theorem 3.1.
\end{Proof}

Let $n,m$ be positive integers, and  $G$ be a graph in $\ell_n^m.$ Clearly, $1\leq m\leq \lfloor\frac{n}{2}\rfloor$.
 Then by Corollary \ref{cor26}, Proposition \ref{prop215} and Theorem \ref{thm31}, we obtain the following result immediately.
\begin{them}\label{thm32}
Let  $G$ be a cactus on $n(\geq 3)$ vertices. Then %$n\geq 3$ be a positive integer and

{\rm (i) } If $n$ is odd,  then $q(G)\leq\frac{(n+2)+\sqrt{n^2-4n+12}}{2}$  with equality if and only if $G \cong H_{\frac{n-1}{2}}^0.$

{\rm (ii) } If $n$ is even,  then $q(G)\leq\frac{(n+1)+\sqrt{n^2-2n+9}}{2}$  with equality if and only if $G \cong H_{\frac{n}{2}-1}^1.$
\end{them}

\end{document}